\documentstyle{article}
\input{amssym.tex}

\newtheorem{Thm}{Theorem}   

\newtheorem{Prop}[Thm]{Proposition}

\begin{document}
\begin{center}
{\Large \bf Decay Of Correlation For  Expanding Toral Endomorphisms\footnote{This paper was pubished in
			 Dynamical Systems,
			  {\em Proceedings of the International Conference in Honor of Professor Liao Shantao
			 Peking University, China, 9 – 12 August 1998.} Ed. L.Wen and Y. P. Jiang, World Scientic,
		1999. However, it is not reviewed in Mathscinet and it is not easy to find the book. So, as asked by colleagues, I put it in ArXiv. Some incomplete references are updated.}
	}\\
\smallskip
Ai-hua FAN \\
\smallskip
Department of Mathematics, University of Picardie, 80039 Amiens,
France \ E-mail: {\tt Ai-Hua.Fan@u-picardie.fr} \\
\medskip
({\it In memory of Professor Liao Shantao})
\end{center}


\begin{abstract}
Let $A$ be an  expanding endomorphism on the torus
${\Bbb T}^d = {\Bbb R}^d /{\Bbb Z}^d$ with its smallest
eigenvalue $\lambda >1$.
Consider the ergodic
system
$({\Bbb T}^d, A, \mu)$ where $\mu$ is Haar measure. We prove
that the correlation  $\rho_{f, g}(n)$ of a pair
of functions $f, g \in L^2(\mu)$ is controlled by the
modulus of $L^2$-continuity $\Omega_{f, 2}(\lambda^{-n})$ and that
 the estimate is  to some extent optimal.
We also prove
the
central limit theorem for the stationary process $f(A^n x)$
defined by a function $f$ satisfying
 $\Sigma_n \Omega_{f,2}(\lambda^{-n}) <\infty$.
An application is given to the Ulam-von Neumann system.

\end{abstract}

\section*{1. Introduction and Main Results}
\setcounter{section}{1} \setcounter{equation}{0}

There has been much interest in the study of  rates of correlation decay
for various kinds of systems \cite{FJ,FP,Liv1,You}
(see the references therein). However few rates for a given class of
test functions are known to be optimal. We intend
in this note to study the simple dynamics of  expanding torus  endomorphisms
and we shall see that optimal rates may be obtained in this case.
Another motivation is to better understand a method introduced
in \cite{FJ,FP} to study
decay of correlations in a general case, which provides rather
precise decay rates but which seems not be able to cover a classical
result that there is an exponential decay rate for H\"{o}lder functions.

We consider the dynamical system $({\Bbb T}^d, A)$
where ${\Bbb T}^d = {\Bbb R}^d /{\Bbb Z}^d$ ($d\geq 1$) is the
$d$-dimensional torus and $A$ is an endomorphism on ${\Bbb T}^d $.
We suppose that $A$ is {\it expanding}, that is,  all of its
eigenvalues have absolute value strictly larger than $1$.
A central theme of the ergodic theory for such a system
is to consider the behavior of $A^n$ as $n \rightarrow \infty$.
A natural way is to describe the behavior of
$A^n$ through an invariant measure.
We take the Haar-Lebesgue  measure  $\mu = dx$ on ${\Bbb T}^d$.
The system $({\Bbb T}^d, A, \mu )$ is strong mixing. That means the
{\it correlation}
$$
    \rho(n) = \rho_{f,g}(n):=
    \int f \cdot g\circ A^n  d\mu - \int f d\mu \cdot \int g d \mu
$$
tends to zero, as
$n \rightarrow \infty$,
 for any $f\in L^2(\mu)$ and any $g \in L^2(\mu)$
(called {\it test functions}).
Our purpose is to study the  rate of correlation decay for a given
pair of test functions $f, g \in L^2(\mu)$.
It is well known that if $g$ and $f$ are H\"{o}lder functions,
the correlation decays exponentially fast. We shall show that for less
regular functions, the correlations decay more slowly and that
different kinds of decay rates are possible.

In order to state our results, we recall here the
modulus of continuity and the modulus of $L^r$-continuity
($1\leq r < \infty$)
of a function $f$ defined on the torus:
$$
    \Omega_f(\delta) = \sup_{|x-y|\leq \delta} |f(x) - f(y)|
$$
$$
    \Omega_{f, r}(\delta) = \sup_{|v|\leq \delta}
                    \|f(\cdot + v) - f(\cdot)\|_r
$$
where $\|f\|_r$ denotes  the norm of $f \in L^r(\mu)$.
Formally, we may write $\Omega_f(\delta)= \Omega_{f, \infty}(\delta)$.

\begin{Thm}
Let $A$ be an expanding integral matrix with
the least eigenvalue (in absolute value) $\lambda>1$.
For $f, g\in L^2(\mu)$, we have
 $$
     |\rho_{f, g}(n)| \leq C \|g\|_2 \Omega_{f, 2}(\lambda^{-n})
     \qquad (\forall n \geq 1)
 $$
where $C>0$ is a constant depending only on $A$.
\end{Thm}

This is actually a consequence of an estimate  on a transfer
operator that we state as the following theorem.
Let $D$ ($\subset {\Bbb Z}^d$)
be a set of representatives of cosets
in ${\Bbb Z}^d/A {\Bbb Z}^d$ (one and only one for each coset).
We call $D$ a set of {\it digits}.  The cardinal of $D$  is equal to
$q = |{\rm det} A|$.
For any  function $f$ on ${\Bbb T}^d$,  define
     $$
     {\mathcal L} f(x) = \frac{1}{q}
          \sum_{\gamma \in D}
           f\left(A^{-1}(x + \gamma ) \right).
     $$
The operator ${\mathcal L}$ is called a {\it transfer operator}.
Note that ${\mathcal L}$ doesn't depend on the choice of
$D$.

\begin{Thm}
Let $A$ be an expanding integral matrix with
the least eigenvalue (in absolute value) $\lambda>1$.
Suppose $\int f d\mu =0$.
If  $f \in L^r(\mu)$ ($1\leq r \leq \infty$), we have
$$
     \|{\mathcal L}^n f
         \|_r
        \leq C \Omega_{f,r}(\lambda^{-n})
        \qquad (\forall n \geq 1).
$$
\end{Thm}

Exponential decays are obtained
for H\"{o}lder continuous functions and functions of  bounded variation
(not necessary continuous). However,
few results are known to be optimal and few work has been carried out for
less regular test functions. As far as we know, an optimal decay is
obtained for some systems studied in \cite{Huh} and less regular test
functions are discussed for expanding systems
in \cite{FJ,FP,Pol}. The method introduced in \cite{FJ,FP}
gives rather precise estimate on the correlations for
functions having Dini continuity ($\int_0^1 \Omega_f(t)/t d t <\infty$).
For the endomorphisms discussed here, it seems that
the modulus of $L^2$-continuity is a good tool to describe
the decay of correlation. It gives a decay rate for {\it every}
pair of test functions in $L^2$ and the decay rate
 is optimal  to some extent, as we shall see by considering
 lacunary trigonometric series.

Another consequence of the above theorem is the following CLT
(Central Limit theorem).

\begin{Thm}
        Let $A$ be an expanding integral matrix and
        let $f\in L^2 ({\Bbb T})$ with $\int f =0$.  Suppose
    $$
          \int_0^1 \frac{\Omega_{f,2}(t)}{t} d t < \infty.
    $$
    Then
    $$
         \frac{1}{\sqrt{n}} \sum_{j=0}^{n-1} f(A^n x)
    $$
    converges in law to a Gaussian variable of zero mean  and
    finite variance
    $$
        \sigma^2 = -\int f^2(x) d x
        + 2 \sum_{n=0}^\infty \int f(x) f(A^n x) d x.
    $$
\end{Thm}

In the one-dimensional  case, M. Kac \cite{Kac}
first proved the CLT for a
function of the class $ {\rm Lip}_\alpha$ with $\alpha >1/2$.
I.A. Ibragimov \cite{Ibr}
weakened the ${\rm Lip}_\alpha$ condition to that the
modulus of $L^r$-continuity $\Omega_{f, r}(\delta)$
of $f$ for some $r >2$ is of order  $\delta^\beta$ ($\beta>0$).
The above Theorem 3 improves  significantly these results.
(But it should be pointed out that a convergence speed for CLT was
obtained in \cite{Ibr}).
For the higher dimensional case,  there were no similar
satisfactory results.

\section*{2. Tiling}

We refer to \cite{GM,LW} for the facts recalled here and for
further information about tilings.
Given a measurable set $T \subset {\Bbb R}^d$, we use
$1_T$ to denote its characteristic function
and $|T|$ to denote its Lebesgue measure.
Given two measurable sets $T$ and $S$,
the notation $T \simeq S$ means that $T$ and $S$ are equal
up to a set of  null Lebesgue measure.

An endomorphism of the torus is represented by an integral matrix.
Suppose $A$ is a $d \times d$ integral matrix which is {expanding},
that is,  all
of its eigenvalues $\lambda_i$ have $|\lambda_i|>1$. Denote
$\lambda = \inf |\lambda_i|$ and
$q = |{\rm det} A|$ ($q\geq 2$ and is an integer).

Take a digit set $D$. Recall that it  consists of
representatives
  of
cosets in ${\Bbb Z}^d/A {\Bbb Z}^d$
(one and only one for a coset).
For each $\gamma \in D$, define
$S_\gamma : {\Bbb R}^d \rightarrow {\Bbb R}^d$
by
     $$
       S_\gamma x = A^{-1} (x + \gamma).
   $$
As for hyperbolic iterated function systems,
it can be  proved that there exists a unique compact set $T$
having the self-affinity
    $$
        T = \bigcup_{\gamma \in D} S_\gamma(T).
    $$
Actually, $|S_{\gamma'}(T) \bigcap S_{\gamma''}(T)|=0$
when $\gamma' \not= \gamma''$.  Therefore the self-affinity
implies
      \begin{equation}
      \sum_{k \in D} 1_T(Ax -k) = 1_T(x)  \qquad {\rm a.e.}
   \end{equation}
It is also known that the compact set $T$ has  the tiling property
   \begin{equation}
      \sum_{k \in {\Bbb Z}^d} 1_T(x -k) = 1 \qquad \qquad {\rm a.e.}
   \end{equation}
Since $T$ satisfies (1) and (2), we say it
generates an {\em integral self-affine tiling}.

The compact set $T$ is called a self-affine tiling {\it tile}.
We have $|T|=1$.
The tiling property allows us to identify ${\Bbb T}^d$
with $T$ up to a null measure set. The self-affinity
allows us to decompose $T$ into $q$ disjoint (up to a
null measure set)  self-affine parts.
If $A$ is a similarity, we say $T$ is  self-similar tiling tile.

\section*{3. Proofs of Theorems}

{\it Notation}: For $\gamma=(\gamma_1, \cdots, \gamma_n) \in D^n$, write
$$
     S_\gamma x = S_{\gamma_n}\circ S_{\gamma_2}\circ
     \cdots \circ S_{\gamma_1}x,
     \qquad
     T_\gamma = S_\gamma(T).
$$
Clearly
$$
     S_\gamma x  = A^{-n}x +A^{-n}\gamma_1 +\cdots + A^{-2}\gamma_{n-1}
          + A^{-1} \gamma_n.
$$
Denoting $b_\gamma =S_\gamma 0$, we get
$$
   {\mathcal L}^n f (x) = \frac{1}{q^n }
                \sum_{\gamma \in D^n}
                f\left( A^{-n}x   + b_\gamma \right).
$$

{\it Proof of Theorem 2}\ \
For $\gamma \in D^n$,
write
$$
     f_\gamma = \frac{1}{|T_\gamma|}\int_{T_\gamma} f(x) d x.
$$
Note that $|T_\gamma| = q^{-n} |T|= q^{-n}$ and that
$$
    \sum_{\gamma \in D^n} f_\gamma =  q^{n} \int_T f
    = q^n \int_{{\Bbb T}^d} f=0.
$$
Then
\begin{eqnarray*}
       {\mathcal L}^n f (x)
       & =   &  \frac{1}{q^n} \sum_{\gamma \in D^n}
               f\left(
                       A^{-n} x + b_\gamma \right)\\
       & =   & \frac{1}{q^n} \sum_{\gamma \in D^n}
               \left[
               f\left(
                       A^{-n} x + b_\gamma \right)
                      - f_\gamma
               \right].
\end{eqnarray*}
Since $T_\gamma = A^{-n}T +b_\gamma$, we have immediately
$$
       |{\mathcal L}^n f (x)|
       \leq  \frac{1}{q^n} \sum_{\gamma \in D^n} \Omega_f({\rm diam}\
       A^{-n}T )
          \leq  C \Omega_f(\lambda^{-n})
$$
where ${\rm diam } B$ denotes the diameter of a set $B$.
We used the fact that ${\rm diam} A^{-n} T \leq a \lambda^{-n}$
for some $a>0$ and the fact that
 $\Omega_f(2\delta) \leq 2 \Omega_f(\delta)$.
The estimate on $\|{\mathcal L}^n f\|_\infty$ is thus proved.

For the estimate on $\|{\mathcal L}^n f\|_r$,
we first use the H\"{o}lder inequality to get
$$
   |{\mathcal L}^n f (x)|^r
       \leq \frac{1}{q^n} \sum_{\gamma \in  D^n}
               \left| f\left(A^{-n} x + b_\gamma \right)
                      - f_\gamma \right|^r  .
$$
By making the change of variables
$ y = A^{-n} x + b_\gamma$,
 we get
$$
  \frac{1}{q^n} \int_T
              \left| f\left(A^{-n} x + b_\gamma \right)
                      - f_\gamma \right|^r  dx
             = \int_{T_\gamma} | f(y) - f_\gamma|^r d y.
$$
Then
$$
   \|{\mathcal L}^n f\|_r^r \leq
          \sum_{\gamma \in D^n}   \int_{T_\gamma} | f(y) - f_\gamma|^r d y.
$$
However
\begin{eqnarray*}
      \int_{T_\gamma} | f(y) - f_\gamma|^r d y
    & =   &   \int_{T_\gamma}
                   \left| \int_{T_\gamma} \left(f(y) - f(x)\right)
                   \frac{dx}{|T_\gamma|}
                  \right|^r d y \\
    &\leq &   \int_{T_\gamma} dy \int_{T_\gamma} |f(y) - f(x)|^r
                 \frac{d x}{|T_\gamma|}\\
    & \leq & \ q^n \int_{T_\gamma} dy \int 1_{ T_\gamma - T_\gamma}(y-x)
                      |f(y) - f(x)|^r dx\\
    & = & \ q^n \int_{T_\gamma} dy \int 1_{ A^{-n} (T-T)}(y-x)
                      |f(y) - f(x)|^r dx
\end{eqnarray*}
Then, if we make a change of variables
$u = y-x$ and $v = x$, we get
\begin{eqnarray*}
     \|{\mathcal L}^n f\|_r^r
    & \leq & \ q^n \int_T dy \int 1_{ A^{-n} (T-T)}(y-x)
                      |f(y) - f(x)|^r dx  \\
    & \leq &  q^n \int_T dv \int_{ A^{-n} (T-T)}
                      |f(u+v) - f(v)|^r du\\
    & = & \ q^n  \int_{ A^{-n} (T-T)} du \int_T
                      |f(u+v) - f(v)|^r dv\\
    & \leq & \ q^n    | A^{-n} (T-T)| \cdot \left[ \Omega_{f, r}({\rm diam}
                      A^{-n}(T - T))
                              \right]^r.
\end{eqnarray*}
Suppose $|T-T| \leq C_1 |T|$ with some sufficiently large
constant $C_1$
(which may depend on the digit set  $D$), we have
$$
 |A^{-n} (T-T)|  \leq C_1 q^n.
 $$
 There is another constant $C_2$ such that
 $$
               {\rm diam} A^{-n}(T-T) \leq C_2 \lambda^{-n}.
 $$
 Therefore for some constant $C_3 >0$, we have
     $$
      \|{\mathcal L}^n f\|_r^r
       \leq C_3
              \left( \Omega_{r, f} (\lambda^{-n}) \right)^r.
     $$
     $\Box$

\bigskip
{\it Proof of Theorem 1}\  \
Assume $\int f =0$, without loss of generality. Consider the operator
$V: L^2 \rightarrow L^2$ defined by
$V f = f \circ A$. It can be verified
that ${\mathcal L}$ is just the adjoint
operator of $V$. Therefore
  $$
 | \rho_{f, g}(n)| = \left|\int g {\mathcal L}^n f \right|
                        \leq \|g\|_2  \|{\mathcal L}^n f\|_2.
  $$
Then Theorem 1 follows from Theorem 2.
$\Box$

    \bigskip
{\it Proof of Theorem 3}\  \
By Theorem 1.1 in \cite{Liv2}, it suffices to verify
$$
    \sum_{n=0}^\infty \left| \int f \cdot f \circ A^n \right|< \infty,
    \qquad
    \sum_{n=0}^\infty  \int |{\mathcal L}^n f | < \infty.
$$
However, both   sums are finite whenever
$\sum_{n=0}^\infty \Omega_{f, 2}(\lambda^{-n}) <\infty$.
This last condition is equivalent to
$\int_0^1 \frac{ \Omega_{f, 2}(s)}{s} ds <\infty$.
$\Box$

\section*{4. Transfer operator, Fourier series and Modulus of continuity}

\begin{Prop}   For $f \in L^1$, we have
$$
    {\mathcal L}^n f(x)
     =   \sum_{k \in {\Bbb Z}^d}
                       \hat{f}(A^{*n} k) e^{2 \pi i
                          \langle k, x\rangle
                                    }
         \qquad (\forall n \geq 1).
$$
\end{Prop}

{\it Proof}\ \
It suffices to prove the expression for $n=1$.
Take a digit set $D^*$ representing ${\Bbb Z}^d /A^*{\Bbb Z}^d $.
Assume that $0 \in D^*$.
Write
$$
   f(x) = \sum_{k \in {\Bbb Z}^d} \hat{f}(k) e^{2 \pi i \langle k, x\rangle}
        = \sum_{k \in {\Bbb Z}^d} \sum_{\beta \in D^*}
            \hat{f}(A^*k+\beta) e^{2 \pi i \langle A^*k+ \beta, x\rangle}.
$$
We have
$$
   f(A^{-1}(x+\gamma))
        = \sum_{k \in {\Bbb Z}^d} \sum_{\beta \in D^*}
            \hat{f}(A^* k+ \beta) e^{2 \pi i
                          \langle A^* k+\beta,\  A^{-1}(x+\gamma)\rangle}
$$
then
\begin{eqnarray*}
    {\mathcal L} f(x)
    & = & \frac{1}{q} \sum_{k \in {\Bbb Z}^d}
                      \sum_{\beta \in D^*}
                      \sum_{\gamma \in D}
                      \hat{f}(A^* k+ \beta) e^{2 \pi i
                          \langle k+ A^{*-1}\beta, \ x+\gamma\rangle
                                    }  \\
    & = & \frac{1}{q} \sum_{k \in {\Bbb Z}^d}
                                e^{2 \pi i \langle k,\ x\rangle}
                      \sum_{\beta \in D^*}
                           \hat{f}(A^* k+ \beta)
                        e^{2 \pi i \langle A^{*-1}\beta,\ x\rangle}
                      \sum_{\gamma \in D}
                       e^{2 \pi i
                          \langle A^{*-1}\beta, \ \gamma\rangle
                                    }
\end{eqnarray*}
So, in order to get
   $$ {\mathcal L} f(x)
     =   \sum_{k \in {\Bbb Z}^d}
                       \hat{f}(A^{*} k) e^{2 \pi i
                          \langle k, x\rangle
                                    }
     $$
 It suffices to note that
$$
    \frac{1}{q}   \sum_{\gamma \in D}
               e^{2 \pi i  \langle A^{*-1}\beta,\ \gamma\rangle}
               = \left\{ \begin{array}{ll}
                              1 & \quad {\rm if} \ \
                                    \beta =0 \\
                              0 & \quad {\rm if} \  \
                                    \beta \not=0.
                        \end{array}
               \right.
$$
In fact, suppose the above sum is not zero.
We have only to show that $\beta =0$.
Since the group $D$
is a product of cyclic groups and
$m^{-1}\sum_{j=0}^{m-1} e^{2 \pi i j x} = 1 $ or $0$ for a real number $x$
according to  $x\in {\Bbb Z}$ or not, we must have
$\langle A^{*-1}\beta,\ \gamma'\rangle$ is an integer for any cyclic
group generator $\gamma'$. Then
$\langle A^{*-1}\beta,\ \gamma\rangle$  is an integer for
any $\gamma \in D$. Let $z \in {\Bbb Z}^d$. Write
$z = \gamma + A k$ with  $\gamma \in D$ and $k \in {\Bbb Z}^d$.
Then
$$
    \langle A^{*-1}\beta,\ z\rangle  =
      \langle A^{*-1}\beta,\ \gamma\rangle
      +  \langle \beta,\ k \rangle = 0    \quad{\rm (mod \ {\Bbb Z})}.
$$
It follows that $\beta = 0$ (mod $A^*{\Bbb Z}^d$). So $\beta = 0$.
$\Box$

\bigskip
{\it Notation}:  " $a_n \approx b_n$" means there are
constants $C_1>0, C_2>0$ such that $C_1 a_n \leq b_n \leq C_n a_n$
for all $n \geq 1$.
\bigskip

From  Theorem 2 and Proposition 1, we get
$$
     \|{\mathcal L}^n f\|_2
     \approx \sqrt{\sum_{k \in {\Bbb Z}^d} |\widehat{f}(A^{*n} k)|^2 }
     \leq C \Omega_{f, 2}(\lambda^{-n}).
$$
$$
      \|{\mathcal L}^n f\|_\infty
      \leq \sum_{k \in {\Bbb Z}^d} |\widehat{f}(A^{*n} k)|
      \leq C \Omega_{f, 2}(\lambda^{-n}).
    $$
These inequalities may become "$\approx$" for some functions $f$.
 Let us consider the  class of functions defined by
 lacunary trigonometric series
 $$
       H (x)= \sum_{k=1}^\infty
       a_n e^{ 2 \pi i
            \langle A^{*k} h, \  x  \rangle}
            \qquad ( h \in {\Bbb Z}^d \setminus \{0\}).
 $$
 Since $\{ A^{*k} h\}_{k\geq 1}$ is a Sidon set \cite{Kah}, we have
 $\|{\mathcal L}^n H\|_r \approx \|{\mathcal L}^n H\|_2 $
 ($\forall 1\leq r <\infty$) and
 \begin{equation}
     \|{\mathcal L}^n H \|_\infty
     \approx \sum_{k=n+1}^\infty |a_{k}| \leq C \Omega_H(\lambda^{-n})
 \end{equation}
  \begin{equation}
         \|{\mathcal L}^n H \|_2
    \approx  \sqrt{
     \sum_{k=n+1}^\infty |a_{k}|^2}
      \leq C \Omega_{H, 2} (\lambda^{-n}).
  \end{equation}
Now let us estimate the modulus of continuity of $H$ by its
Fourier coefficients.
The following proposition is immediate, just because
$$
     |H(x+ y) -H(x)|
     \leq \sum_{k=1}^n
       |a_k| \left| e^{ 2 \pi i
            \langle A^{*k} h, \  y  \rangle} - 1 \right|
     + \sum_{k=n+1}^\infty
       |a_k|
$$
$$
     \int |H(x+ y) -H(x)|^2 d x
     \leq \sum_{k=1}^n
       |a_k|^2 \left| e^{ 2 \pi i
            \langle A^{*k} h, \  y  \rangle} - 1 \right|^2
     + \sum_{k=n+1}^\infty
       |a_k|^2.
$$

 \begin{Prop} Let $A$ be an expanding integral similarity
 matrix with spectral radius $\lambda >1$.
 Let $H$ be the function defined by the above lacunary
 trigonometric series. We have
 $$
        \Omega_H(\lambda^{-n})
        \leq  C
        \lambda^{-n}\sum_{k=1}^n  |a_k|\lambda^k +
        \sum_{k=n+1}^\infty  |a_k|
 $$
 $$
        \Omega_{H,2}(\lambda^{-n})^2
        \leq
          C \lambda^{-2 n}\sum_{k=1}^n  |a_k|^2\lambda^{2k} +
              \sum_{k=n+1}^\infty  |a_k|^2
 $$
  where $C>0$ is a constant.
 \end{Prop}

 If  $a_n = \frac{1}{n^\alpha}$ with $\alpha >1$, then
 $$
     \|{\mathcal L}^n H\|_\infty
     \approx  \Omega_H (\lambda^{-n})  \approx \frac{1}{n^{\alpha -1}}
$$
$$
         \|{\mathcal L}^n H\|_2
      \approx \Omega_{H, 2} (\lambda^{-n})  \approx \frac{1}{n^{\alpha -1/2}}.
 $$

 If  $a_n = \frac{1}{n \log^\beta n}$ with $\beta >1$, then
 $$
       \|{\mathcal L}^n H\|_\infty
      \approx \Omega_H (\lambda^{-n})
      \approx \frac{1}{(\log n)^{\beta -1}}
$$
$$
       \|{\mathcal L}^n H\|_2
      \approx
      \Omega_{H, 2} (\lambda^{-n})
       \approx \frac{1}{ (\log n)^{\beta-1/2}}.
 $$

 If  $a_n = \theta^n$ with $\frac{1}{\lambda}< \theta <1$, then
 $$
        \|{\mathcal L}^n H\|_\infty
        \approx
        \|{\mathcal L}^n H\|_2
      \approx \Omega_H(\lambda^{-n})
      \approx \Omega_{H,2}(\lambda^{-n}) \approx \theta^n.
 $$

In fact, we have only to use (3) and (4) for getting lower bounds
and to use Proposition 2 for getting upper bounds.
When $a_n = \frac{1}{n^\alpha}$, it suffices to remark that
(for any $\lambda >1$ and any $\alpha >1$)
$$
\sum_{k=1}^n  \frac{\lambda^k}{k^\alpha} \approx
           \sum_{k=1}^{[n/2]}        \frac{\lambda^k}{k^\alpha}
         + \sum_{k=[n/2] +1}^n       \frac{\lambda^k}{k^\alpha}
           \leq C\left( \lambda^{n/2} + \frac{\lambda^n}{n^\alpha}\right)
$$
$$
   \sum_{k=n+1}^\infty  \frac{1}{k^\alpha} \approx \frac{1}{n^{\alpha-1}}
$$
where $[n/2]$ denotes the integral part of $n/2$.
In the same way, we treat the case $a_n = \frac{1}{n \log^\beta n}$.
The case $a_n = \theta^n$ is simpler.
More generally, suppose $|a_n|$ is a decreasing sequence such that
$$
    \limsup_{n \rightarrow \infty} \frac{|a_{[\delta n]}|}{|a_n|}< \infty,
\qquad
    \limsup_{n \rightarrow \infty} \frac{\lambda^{-(1-\delta)n}}{|a_n|}< \infty
$$
for some $0<\delta<1$. Then we have
$$
   \Omega_H(\lambda^{-n}) \approx \sum_{k= n+1}^\infty |a_k|,
\qquad
   \Omega_{H,2}(\lambda^{-n}) \approx \sqrt{ \sum_{k= n+1}^\infty |a_k|^2 }.
$$

\bigskip
Let us finish our discussion by making some remarks:\\

        1.  In all  three cases of $H$  discussed above,
        the estimate provided by Theorem 2 is optimal.
        Let us point out  that
        $H$
        belongs to the class of functions with
        $\Omega_f(\delta) =O(1/|\log \delta|^{\alpha-1})$
        when $a_n= \frac{1}{n^\alpha}$;
        $H$
        belongs to the class of functions with
        $\Omega_f(\delta) =O(1/(\log |\log \delta|)^{\beta-1} ))$
        when $a_n= \frac{1}{n \log^\beta n}$;
        $H$
        belongs to the class of functions with
        $\Omega_f(\delta) =O(\delta^{ \frac{\log \theta}{\log \lambda} })$
        when $a_n= \theta^n$.

        2. For every $f \in L^r$ ($1\leq r <\infty)$,
         $\|f(\cdot +y) -f(\cdot)\|_r$ is continuous
         as a function of $y$. It follows that
         $\lim_{\delta \rightarrow 0} \Omega_{f, r}(\delta)=0$.
         Then, by Theorem 2,
         $\|{\mathcal L}^n f - \int f \|_r$
         tends to zero for any $f \in L^r$.

	3. Let $0<\delta_n <1$ be an arbitrary decreasing sequence
		tending to zero. There is a function $H$ such that
		$\|{\mathcal L}^n H\|_\infty  \approx \delta_n$. In fact,
		 it suffices to take the function $H$ defined by
		$a_n = \delta_{n-1}-\delta_n$ (with $\delta_0=1$).
		Also, if we take the  function $H$ defined by
                $a_n = \sqrt{ \delta_{n-1}^2 - \delta_n^2 }$,
		we have $\|{\mathcal L}^n H\|_2  \approx \delta_n$.

        4. The function $H$ defined above
        with $a_n = \frac{1}{n^\alpha}$ ($1<\alpha\leq 2$)
        or $a_n = \frac{1}{n \log^\beta  n}$ ($\beta >1$)
        is a continuous function but not of summable variation.
        However, $\|{\mathcal L}^n f \|_\infty$ tends to zero
        with a precisely known convergence speed.
        Such a situation was not seen before.

\section*{5. Ulam-von Neumann map}
The map $Uy = 1 -2y^2$ from $I= [-1, 1]$ into itself was studied
by Ulam and von Neumann \cite{UV}. This Ulam-von Neumann map $U$
is conjugate to the tent map $Tx = 1 - 2|x|$. More precisely,
we have $U \circ h = h \circ T$ where $h$ is the conjugacy defined
by $h(x) =  \sin \frac{\pi}{2} x$.   The Lebesgue measure $dx$
(normalized so that $I$ has measure $1$)
is $T$-invariant and its image under $h$,
$d \mu = \frac{2}{\pi} \frac{dy}{\sqrt{1-y^2}}$, is $U$-invariant.
Consider now the system $(I, U, \mu)$. For this system, the transfer
operator is defined by
$$
   {\mathcal L}_U f(y) = \sum_{z \in U^{-1}y} \frac{f(z)}{|U'(z)|}.
$$

\begin{Thm} For any  $f \in {\rm Lip}_\alpha$
($0<\alpha<1$) with $\int f d\mu =0$,
we have $\|{\mathcal L}_U^n f\|_2 \leq C 2^{-\alpha n}$.
For $g(y) = \log |y| + \log 2$, we have  $\int g d \mu =0$ and
$\|{\mathcal L}_U^n g\|_2 \approx  2^{- n}$.
\end{Thm}

{\it Proof}\ \
First observe that Theorem 2 remains true for the system
$(I, T, dx)$ because $I$ can  be decomposed into
$I= S_0(I) \bigcup S_1(I)$,
where $S_0 x = \frac{x-1}{2}, S_1 x = \frac{1-x}{2}$
are inverses of $T$.
Let ${\mathcal L}_T$ be the transfer operator
associated to $T$, which is defined in the same way as  ${\mathcal L}_U$.
By using the fact that ${\mathcal L}_U$
is the adjoint operator of $f \rightarrow f\circ U$ acting
on $L^2(\mu)$ and the similar fact about ${\mathcal L}_T$,
we get the relation between ${\mathcal L}_U$
and ${\mathcal L}_T$:  for any $g, f \in L^2(\mu)$
\begin{eqnarray*}
   &  &  \int g \cdot {\mathcal L}_U^n f d \mu
     =
     \int g \circ U^n \cdot  f d \mu \\
   &  =  & \int g \circ U^n \circ h \cdot  f \circ h dx
     =  \int g  \circ h \circ T^n  \cdot  f \circ h dx  \\
   &  = &  \int g  \circ h   \cdot   {\mathcal L}_T^n (f \circ h) dx
     =   \int g     \cdot   {\mathcal L}_T^n (f \circ h)\circ h^{-1} d\mu.
\end{eqnarray*}
It follows that
$\|{\mathcal L}_U^n f\|_2 = \|{\mathcal L}_T^n (f\circ h)\|_2$. In fact,
\begin{eqnarray*}
  & &  \int {\mathcal L}_U^n f \cdot \overline{ {\mathcal L}_U^n f  } d \mu
   = \int {\mathcal L}_U^n f \cdot
            {\mathcal L}_T^n (\overline{f}\circ h)\circ h^{-1} d \mu \\
  & =  &\int {\mathcal L}_T^n (f \circ h) \circ h^{-1} \cdot
            {\mathcal L}_T^n (\overline{f}\circ h)\circ h^{-1} d \mu \\
  &  = & \int {\mathcal L}_T^n (f \circ h)  \cdot
            {\mathcal L}_T^n (\overline{f}\circ h) d x.
\end{eqnarray*}

So, we have only to estimate $\|{\mathcal L}_T^n(f \circ h)\|_2$.
To this end, apply Theorem 2 ( its variant mentioned above).
We  are then led to estimate the modulus of continuity of
$f\circ h$.  Since $h$ is Lipschitz, we have
$\Omega_{f\circ h, 2}(\delta)\leq \Omega_{f}(C' \delta)\leq C
\Omega_{f}( \delta)$.

However, $\log |y|$ is neither continuous nor of bounded variation.
But it is in $L^2(\mu)$,
equivalently $\log |h(x)| \in L^2(dx)$. Note that
$$
 \log |h(x)| =\log |\sin \frac{\pi}{2} x|
   = - \log 2 - \sum_{n=1}^\infty \frac{\cos \pi n x}{n}.
$$
It follows that $\int \log |y| d \mu = \int \log |\sin \frac{\pi}{2} x| dx
= -\log 2$.
From the above series and Proposition 1,
we get $\|{\mathcal L}^n_U g\|_2 = \|{\mathcal L}^n_T g\circ h\|_2
\approx 2^{-n}$.
$\Box$

\bigskip
Remark that we can also get
$\Omega_{g\circ h, 2}(\delta)\leq C \sqrt{\delta}$.
In fact, for any $u \not = 0$,  take $N$ to be the integral part of $1/|u|$.
According to Parseval's equality,
 we have
\begin{eqnarray*}
  & & \int [ \log |h (x+u)|-
      \log |h(x)| ]^2 d x \\
  &= &2 \sum_{n=1}^\infty  \frac{\sin^2 \frac{\pi}{2} n u}{n^2}
 \leq  \frac{\pi^2}{2}N |u|^2 +   2 \sum_{n=N +1}^\infty \frac{1}{n^2}
      \leq C \sqrt{|u|}.
\end{eqnarray*}

\bigskip

It is known that the Liapunov exponent of $(I, U, \mu)$ is  equal to
$$
\int \log |U'(y)| d\mu(y) = \log 2.
$$That means for $\mu$
almost all points $y$,  $\frac{1}{n} \log |(U^n)'(y)|$
converges to $\log 2$. As a consequence of the last theorem,
we get that $\frac{1}{\sqrt{n}} [\log |(U^n)'(y)| - n \log 2]$
converges in law to a centered gaussian variable
(following the arguments in the proof of Theorem 3).

\bigskip

The above discussion on Ulam-von Neumann map reveals a possibility
to reduce the study of a general system to that of an endomorphism
on torus or a system like the tent map.
It is the case  when there is a conjugacy
between the two systems and when the
conjugacy has some smoothness so that the modulus of continuity
of $f \circ h$ is small.

\bigskip {\it Acknowledgement} \ \ The author would like to thank
Yunping JIANG for his valuable discussions and Oliver JENKINSON
for his careful reading of an earlier version of the paper.

\end{document}